\documentclass[11pt,a4paper,oneside,article]{article}
\usepackage[T2A]{fontenc} \usepackage[utf8]{inputenc} \usepackage{amsfonts,
amsmath, amsthm} \usepackage{amssymb} \usepackage{algorithm}
\usepackage{algorithmic} \usepackage{xcolor} \usepackage{url}
\usepackage{etoolbox} \usepackage{pdfpages}

\newcommand*{\norm}[1]{\left\lVert#1\right\rVert}
\newcommand{\argdualmin}{\mathop{\mathrm{argdual}}}
\DeclareMathOperator*{\Argmin}{Arg\,min}

\newtheorem{definition}{Definition} \newtheorem{theorem}{Theorem}
\newtheorem{assumption}{Assumption}

\def\E{\mathbb E}

\providetoggle{thesis} \settoggle{thesis}{false}

\begin{document}

\begin{titlepage}
   \begin{center}
National Research University Higher School of Economics
       \vspace{0.5cm}
       \begin{flushright}
       \it
as a manuscript
       \end{flushright}
       
       \vspace{1.5cm}
       
       Alexander Igorevich Tyurin
       
\vspace{1.5cm}       
\textbf{\Large Development of a method for solving structural optimization problems}
       
       \vspace{4.5cm}
       
       PhD Dissertation Summary
       
       for the purpose of obtaining academic degree
       
       Doctor of Philosophy in Computer Science

       \vfill
            
       Moscow - 2020
            
   \end{center}
\end{titlepage}
The PhD dissertation was prepared at International Laboratory of Stochastic Algorithms and High-Dimen\-sional Inference, National Research University Higher School of
Economics \\

{\bf Academic Supervisor:} \\

Alexander Vladimirovich Gasnikov, Doctor of Sciences in Mathematical Modelling, Numerical Methods and Software Complexes, Senior Research Fellow at International Laboratory of Stochastic Algorithms and High-Dimen\-sional Inference, National Research University Higher School of
Economics

\newpage
\section{Introduction}

Optimization methods have a significant impact on all spheres of human society. It is difficult to list all recent activities where optimization methods are used to solve practical problems. In many problems of economics, engineering, programming optimization methods are helpful. Optimization methods came up with computer engineering in the twentieth century. That is when the active development of the modern theory of optimization began. The pioneer is L. Kantorovich \cite{kantarovich1939, polyakhistory}, who considered linear programming problems in engineering and economics. In the 50s-60s, cutting edge works were done by G. Rubinstein, E. Ventsel, N. Vorobyov, D. Yudin, E. Golstein, N. Shor, B. Polyak. Yu. Ermoliev, L. Pontryagin, etc. In those years, researchers proposed the following methods: Pontryagin's maximum principle, projection gradient method, cutting plane method, penalty method, Newton's method for constrained optimization, subgradient method, the center of gravity method, etc. In 70th, A. Nemirovskii and D. Yudin have a significant impact on optimization development with work \cite{nemirovkii1979}. In this work, they used the oracle concept (black box), which for any input point returns, for example, function and gradient value. A. Nemirovski and D. Yudin obtained lower bounds for convergence rates for some general optimization problems classes (convex optimization problems, optimization problems with Lipchitz continuous functions, smooth strongly convex optimization problems, etc.). We should note that optimization methods that achieve corresponding lower bounds are proposed later for some classes of problems. In particular, Yu. Nesterov developed the fast gradient method \cite{nesterov1983}, which has the convergence rate inversely proportional to the root of the accuracy of the solution by function in the class of functions with Lipchitz continuous gradient. This convergence rate is optimal in the sense of black-box oracle calls. The same result was obtained for smooth strongly convex functions.

It means that for many classes of optimization problems, optimal methods were developed; however, the progress did not stop. In practice, optimization tasks have some structure that allows developing new algorithms for every problem with faster convergence rates. Let us note some popular examples from structural optimization. Composite optimization solves optimization problems that can be represented as a sum of smooth and nonsmooth functions. Despite the fact that the sum is a nonsmooth function, with some additional assumptions about the nonsmooth part, we can develop methods that have convergence rates as in smooth optimization tasks \cite{nesterov2013}. Similarly, using the structure of optimization tasks, we can propose algorithms with more optimistic convergence rates for the following optimization problems: functions with Holder continuous gradients \cite{nesterov2015universal}, superposition of functions (min-max problems) \cite{nemirovski1995information, lan2015bundle}, transportation problems \cite{baymurzina2019universal, gasnikov2018dualtransport, gasnikoveffectivnie}, clustering by electorial model \cite{nesterov2018soft}, etc.

We should note another development direction in structural optimization theory that connected with different requirements about an oracle. In general, an oracle is a black-box framework that makes calculations. The complexity of optimization methods is estimated by the number of calls of an oracle. In classical optimization theory \cite{nemirovkii1979, nesterov2004} oracles for some query point return a function value (zero-order oracle), gradient / subgradient (first-order oracle), hessian (second-order oracle), etc. In particular, for a smooth convex function $f$ with $L$--Lipchits gradient it is known that there is exist an oracle \cite{nesterov2004} such that for some query point $y$ returns a pair $(f(y), \nabla f(y))$ and the following inequality holds:
\begin{equation}
\label{introduction:smooth_oracle_inequality}
    0 \leq f(x) - (f(y) + \langle \nabla f(y), x - y\rangle) \leq \frac{L}{2}\norm{x - y}^2 \quad \forall x \in Q,
\end{equation}
where $Q$ is a convex closed set on which a function $f$ is defined. We have the left inequality from the convexity of a function $f$ and the right inequality we have from $L$--Lipchits gradient assumption. Using this oracle, we can obtain the optimal convergence rate for this class of optimization problems. In detail, after $N$ calls of the oracle, the guaranteed convergence rate is equal to $\mathcal{O}\left(\frac{L R^2}{N^2}\right)$ \cite{nesterov1983, nesterov2004}, where a constant $R$ is a distance between a method's starting point and the closest optimal point. In practice and theory, inequalities \eqref{introduction:smooth_oracle_inequality} do not always hold. Even for smooth optimization problems with Lipchits gradients, we can not get precise values of function and gradient due to calculation errors or the fact, that we obtain these values using some auxiliary problem. For these examples, inequalities \eqref{introduction:smooth_oracle_inequality} do not hold. Indeed, we can show that \eqref{introduction:smooth_oracle_inequality} holds only for some unique pair $(f(y), \nabla f(y))$. In \cite{devolder2011} authors introduced $(\delta, L)$--oracle that for some query point $y$ returns a pair $(f_\delta(y), \nabla f_\delta(y))$ such that
\begin{equation}
\label{introduction:delta_l_oracle}
    0 \leq f(x) - (f_\delta(y) + \langle \nabla f_\delta(y), x - y\rangle) \leq  \frac{L}{2}\norm{x - y}^2 + \delta \quad \forall x \in Q.
\end{equation}
Unlike \eqref{introduction:smooth_oracle_inequality}, a pair $(f_\delta(y), \nabla f_\delta(y))$ is not unique and, in general, is not equal to $(f(y), \nabla f(y))$. The proposed oracle allows us to generalize the classical gradient and fast gradient type methods to a wider class of tasks. From \cite{devolder2011}, the guaranteed convergence rate for the fast gradient method is equal to $\mathcal{O}\left(\frac{L R^2}{N^2} + N \delta \right)$ and for the gradient method is equal to $\mathcal{O}\left(\frac{L R^2}{N} + \delta \right)$ using $(\delta, L)$--oracle. We should note, that obtained convergence rates do not require smoothness of optimization problems. We can obtain $(\delta, L)$--oracle for the following optimization problems \cite{devolder2011}: functions with Holder continuous subgradient, functions obtained by smoothing techniques \cite{nesterov2005, nesterov2005excessive, nesterov2007smoothing}, Moreau-Yosida regularization \cite{lemarechal1997}, and composite optimization \cite{nesterov2013}. Note that there is another concept of inexact oracle \cite{aspremont} that is a particular case of $(\delta, L)$--oracle \cite{devolder2011}.

A valuable property of optimization methods is to effectively resolve a dual solution from a primal solution (or vice versa)  \cite{anikin2017dual, boyd2004convex, nesterov2015complexity, nesterov2009primal}, which is called primal--duality. This property is advantageous in transportation problems \cite{baymurzina2019universal, gasnikov2018dualtransport, gasnikoveffectivnie}, machine learning problems (ex. SVM), etc. Another useful property of optimization methods is to be robust to cases when instead of a gradient oracle returns a stochastic gradient, random, unbiased vector w.r.t. a real gradient. Stochastic optimization methods are trendy because they allow reducing the calculation cost of a descent direction. It is impossible for some optimization problems \cite{goodfellow2016deep, krizhevsky2012imagenet} to calculate a gradient in a reasonable amount of time, even at one point.
Therefore, in this thesis, we focus on the extension of our results to primal--duality and stochasticity.

{\bf  Object and goals of the dissertation:} unification of previously proposed gradient-type methods into one method using a special concept of inexact model. 
% Основным требованием к неточной модели функции является возможность использования ее для как можно большего множества оптимизационных задач. 
Develop a series of methods that can solve generalized optimization problem statements and use its structure with the aid of the proposed concept of inexact model. Moreover, prove corresponding rates of convergence, when possible, in an optimal way for some classes of optimization problems. 

{\bf  The obtained results:}
\begin{enumerate}
    \item We propose concepts of inexact model for gradient-type methods.
    \item We developed the adaptive gradient and fast gradient methods for optimizations tasks that support the concept of inexact model ($(\delta, L,\\ \norm{})$--model).
    \item We developed the gradient method for optimization tasks with relative smoothness that support the concept of inexact model ($(\delta, L,\\ V)$--model).
    \item We developed the primal--dual adaptive gradient and fast gradient methods for optimizations tasks that support the concept of inexact model ($(\delta, L, \norm{})$--model).
    \item We developed the stochastic gradient and fast gradient methods for stochastic optimizations tasks that support the concept of inexact model ($(\delta_1, \delta_2, L, \norm{})$--model).
    \item We propose the heuristic (without theoretical guarantees) adaptive stochastic fast gradient method that is based on the adaptive fast gradient method and the stochastic fast gradient method.
\end{enumerate}

{\bf Author’s contribution} includes the development of optimization methods for oracles that return an inexact model of a function, proving convergence rates of corresponding methods, and constructing of inexact models for problems from structural optimization. The heuristic adaptive algorithm is proposed for stochastic optimization problems.

{\bf Novelties:} we developed the adaptive gradient and fast gradient methods for oracles that return an inexact model of a function. Further, we constructed the gradient method for problems with relative smoothness, the primal--dual adaptive gradient and fast gradient methods, and the stochastic nonadaptive gradient methods that support an inexact model of a function. We attempted to enrich the stochastic nonadaptive fast gradient method with adaptivity. However, we were only able to develop the heuristic adaptive fast gradient method that showed high performance in practice.

As a result of the work of this thesis, 8 papers were published:

{\bf First-tier publications:}

\begin{enumerate}
    \item Gasnikov A., Tyurin A. (2019) Fast gradient descent for convex mi\-ni\-mi\-zation problems with an oracle producing a $(\delta, L)$-model of function at the requested point. Computational Mathematics and Mathematical Physics, 59, 7, 1085--1097, Scopus Q2 (main co-author; the author of this thesis formulated and proved convergence rate theorems for the gradient and fast gradient methods (Theorem 1 and 2), presented a description of examples (Section 4)).
    \item Stonyakin F., Dvinskikh D., Dvurechensky P., Kroshnin A., Kuznetsova O., Agafonov A., Gasnikov A., Tyurin A., Uribe С., Pasechnyuk D., Artamonov S. (2019) Gradient methods for problems with inexact model of the objective. Lecture Notes in Computer Science, 11548, 97--114, Scopus Q2 (the author of this thesis prepared the text of section 2 and proved the convergence rate theorem for the gradient method for optimization problems with relative smoothness (Theorem 1)).
    \item Ogaltsov A., Tyurin A. (2020) A heuristic adaptive fast gradient method in stochastic optimization problems. Computational Mathematics and Ma\-the\-matical Physics, 60, 7, 1108--1115, Scopus Q2 (main co-author; the author of this thesis proposed the heuristic adaptive fast gradient method for stochastic optimization problems (Algorithm 2), did an analysis and justification).
    \item Dvurechensky P., Gasnikov A., Omelchenko A., Tyurin A. A stable alternative to Sinkhorn’s algorithm for regularized optimal transport. Lecture Notes in Computer Science, 12095, 406--423, Scopus Q2 (the author of this thesis helped with the development of Algorithm 1 and the proof of Theorem 1).
    \item Dvinskikh D., Omelchenko A., Gasnikov A., Tyurin A. Accelerated gradient sliding for minimizing the sum of functions. Doklady Mathematics, 101, 3, Scopus Q2 (in press), (the author of this thesis helped with the text of this paper and the proofs of intermediate results).
\end{enumerate}

{\bf Second-tier publications:}

\begin{enumerate}
   \item Tyurin A. (2020) Primal-dual fast gradient method with a model. Computer Research and Modeling, 12, 2, 263--274, Scopus Q3.
    \item  Dvinskikh D., Tyurin A., Gasnikov A., Omelchenko S. (2020) Accelerated and nonaccelerated stochastic gradient descent with model conception. Mathematical Notes, 108, 4, Scopus Q3 (in press), (main co-author; the author of this thesis developed the fast gradient method for stochastic optimization tasks, provided description of examples).
    \item Anikin A.,  Gasnikov A., Dvurechensky P., Tyurin A., Chernov A. (2017) Dual approaches to the minimization of strongly convex functionals with a simple structure under affine constraints. Computational Mathematics and Mathematical Physics, 57, 8, 1262--1276, Scopus Q3. (the author of this thesis helped in writing of remarks).
\end{enumerate}

{\bf Reports at conferences and seminars:}

\begin{enumerate}
    \item 8th Moscow International Conference on Operations Research, Russia, Moscow. (17.10.2016 - 22.10.2016). Dual fast gradient method for entropy--linear programming problems.
    \item 59th MIPT Scientific Conference, Russia, Dolgoprudny. (21.11.2016 - 26.11.2016). Adaptive fast gradient method for convex min--max problems.
    \item Workshop ``Three oracles'', Russia, Skolkovo. (28.12.2016).
    On several extensions of similar triangles method.
    \item Scientific conference "Modeling the Co-evolution of Nature and Society: problems and experience" devoted to the 100-th anniversary of N. N. Moiseev, Russia, Moscow. (7.11.2017 - 10.11.2017). Adaptive similar triangles method and its application in calculation of regularized
optimal transport.
    \item 60th MIPT Scientific Conference, Russia, Dolgoprudny. (20.11.2017 - 25.11.2017). The mirror triangle method with a generalized inexact oracle.
    \item The 23rd International Symposium on Mathematical Programming, France, Bordeaux. (1.7.2018 - 6.7.2018).  Universal Nesterov’s gradient method in general model conception.
    \item 62th MIPT Scientific Conference, Russia, Dolgoprudny. (18.11.2019 - 23.11.2019). Primal--dual fast gradient method with a model.
\end{enumerate}

The reported study was funded by RFBR, project number 19-31-90062 and project number 18-31-20005 mol-a-ved.

\section{Convex optimization problem}
Let us describe the mathematical formulation of a convex optimization problem \cite{nesterov2004}. Given an objective function $f(x) : Q \longrightarrow \mathbb{R}$, a set $Q$ is a subset of finite--dimensional linear vector space $\mathbb{R}^n$, and a norm $\norm{\cdot}$ in $\mathbb{R}^n$. Conjugate norm we define as $$\norm{\lambda}_* = \max_{\norm{\nu} \leq 1, \nu \in \mathbb{R}^n} \langle \lambda, \nu \rangle \quad \forall \lambda \in \mathbb{R}^n.$$

Let us define prox--function and Bregman divergence \cite{bregman1967relaxation}, \cite{ben-tal2015lectures} (p. 327):

\begin{definition}
$d(x):Q \rightarrow \mathbb{R}$ is a prox--function if $d(x)$ is continuously differentiable on $\textnormal{int }Q$ and a function $d(x)$ is 1-strongly convex w.r.t. a norm $\norm{}$ on $\textnormal{int }Q$.
\end{definition}
\begin{definition}
A function
\begin{align*}
V[y](x) = d(x) - d(y) - \langle\nabla d(y), x - y\rangle
\end{align*}
is called Bregman divergence, where $d(x)$ is a prox--function.\footnote{In paper \cite{gasnikov2019fast}, $V[y](x)$ is denoted as $V(x, y)$.}
\end{definition}
The introduction of Bregman divergence allows us to obtain more general results for convergence rates. The classical example of Bregman divergence is the function $V[y](x) = \frac{1}{2}\norm{x - y}_2^2$.

Further, we assume that
\begin{enumerate}
  \item A set $Q \subseteq \mathbb{R}^n$ is a convex and closed set.
  \item A function $f(x)$ is continuous and convex on $Q$.
  \item A function $f(x)$ is lower bounded on $Q$ and attains its minimum at some point from $Q$ (not necessarily unique).
\end{enumerate}
We consider the following optimization problem:
\begin{align}
\label{mainTask}
f(x) \rightarrow \min_{x \in Q}.
\end{align}
A point $x_*$ is a solution of the optimization problem if inequality
$f(x_*) \leq f(x)$ holds for all $x \in Q$. Also, we call a point $x_{\varepsilon}$ as $\varepsilon$--solution if $f(x_{\varepsilon}) - f(x_*) \leq \varepsilon$ for all $x \in Q$. The main task of numerical convex optimization is to find $\varepsilon$--solution.

\subsection{The concept of inexact solution}

Now we define the concept of inexact solution (see \cite{ben-tal2015lectures}) that we use in our methods.
\begin{definition}
\label{solNemirovskiy}
Given an optimization problem
\begin{gather*}
\psi(x) \rightarrow \min_{x \in Q},
\end{gather*}
where $\psi(x)$ is a convex function, then we denote $\Argmin_{x \in Q}^{\widetilde{\delta}}\psi(x)$ by a set of $\widetilde{x}$ such that
\begin{gather*}
\exists h \in \partial\psi(\widetilde{x}), \,\,\, \langle h, x - \widetilde{x}  \rangle \geq -\widetilde{\delta} \,\,\,\, \forall x \in Q,
\end{gather*}
where $\partial\psi(\widetilde{x})$ is a subderivative of a function $\psi$ at a point $\widetilde{x}$.
Any point from $\Argmin_{x \in Q}^{\widetilde{\delta}}\psi(x)$ we denote as $\arg\min_{x \in Q}^{\widetilde{\delta}}\psi(x)$.

\end{definition}

This definition is stringent compared to the definition of $\delta$--solution by function (see \cite{gasnikov2019fast}). For instance, both definitions are equivalent when $\widetilde{\delta} = 0$. However, in some more general cases, we can derive solutions in terms of Definition \ref{solNemirovskiy} from $\delta$--solution by function (see \cite{stonyakin2020inexact, gasnikov2019fast}).

\section{Contents}

In this section, we describe results and conclusions in more detail.

\subsection{Inexact model of a function}
In \cite{devolder2011}, the authors proposed $(\delta, L)$--oracle and corresponding methods that allow solving a vast number of optimization problems. Let us introduce concepts of inexact model ($(\delta, L)$--model) of a function that generalizes $(\delta, L)$--oracle.
\begin{definition}
\label{definition:model}
  Given a convex function $\psi_{\delta}(x, y)$ w.r.t. $x$ on a set $Q$ such that $\psi_{\delta}(x, x) = 0$ for all $x \in Q$. The function $\psi_{\delta}(x, y)$ is $(\delta, L, \norm{})$--model of a function $f$ at a point $y$ w.r.t. $\norm{\cdot}$ with value $f_\delta(y)$ if
  for all $x \in Q$ inequalities
  \begin{gather}
  \label{model_def}
  0 \le f(x) - (f_{\delta}(y) + \psi_{\delta}(x, y)) \le \frac{L}{2}\norm{x - y}^2 + \delta
  \end{gather}
  hold for some values $L, \delta \geq 0$.\footnote{In papers \cite{gasnikov2019fast, tyurin2020primal}, $(\delta, L, \norm{})$--model is defined as $(\delta, L)$--model.}
\end{definition}

From the view of an oracle concept, we can assume that for a query point $y$, the oracle returns a pair $(f_\delta(y), \psi_{\delta}(x, y))$.
Also, we can provide a more general definition by using so--called relative smoothness \cite{bauschke2016, lu2018, stonyakin2019gradient}:
\begin{definition}
\label{definition:model_relative}
  Given a convex function $\psi_{\delta}(x, y)$ w.r.t. $x$ on a set $Q$ such that $\psi_{\delta}(x, x) = 0$ for all $x \in Q$. The function $\psi_{\delta}(x, y)$ is $(\delta, L, \norm{})$--model of a function $f$ at a point $y$ w.r.t. Bregman divergence $V$ with value $f_\delta(y)$ if
  for all $x \in Q$ inequalities
  \begin{gather}
  \label{model_def_relative}
  0 \le f(x) - (f_{\delta}(y) + \psi_{\delta}(x, y)) \le LV[y](x) + \delta
  \end{gather}
  hold for some values $L, \delta \geq 0$.\footnote{In paper \cite{stonyakin2019gradient} $(\delta, L, V)$--model is defined as $(\delta, L)$--model.}
\end{definition}

We can obtain Definition \ref{definition:model} from Definition \ref{definition:model_relative} if we take Bregman divergence $V[y](x) = \frac{1}{2}\norm{x - y}^2$.
The oracle that produces $(\delta, L, \norm{})$--model from Definition \ref{definition:model} is more universal than $(\delta, L)$--oracle (see \eqref{introduction:delta_l_oracle}). Indeed, it is enough to take $\psi_{\delta}(x, y) = \langle \nabla f_\delta(y), x - y \rangle$.

In \cite{dvinskikh2020accelerated}, we propose a more general definition than Definition \ref{definition:model} by the introduction of additional noise, namely:
\begin{definition}
\label{definition:model_stochastic}
  Given a convex function $\psi_{\delta}(x, y)$ w.r.t. $x$ on a set $Q$ such that $\psi_{\delta}(x, x) = 0$ for all $x \in Q$. The function $\psi_{\delta}(x, y)$ is $(\delta_1, \delta_2, L, \norm{})$--model of a function $f$ at a point $y$ w.r.t. $\norm{\cdot}$ if
  for all $x \in Q$ inequalities
  \begin{gather}
  \label{model_def_stochastic}
  -\delta_1(x,y) \le f(x) - (f(y) + \psi_{\delta}(x, y)) \le \frac{L}{2}\norm{x - y}^2 + \delta_2
  \end{gather}
  hold for some values $L \geq 0$, $\delta_2$, and $\delta_1(x,y)$.\end{definition}
We can show, that $(\delta, L, \norm{})$--model is $(\delta, \delta, L, \norm{})$--model with $\delta_1(x,y) = \delta$ and $\delta_2 = \delta$. This concept is helpful in stochastic optimization problems. (see Section \ref{subsection:stochastic_gradient_descent}).

\subsection{Examples of inexact models of a function}
Let us provide some examples of inexact models for different optimization tasks.
\begin{enumerate}
\item {\it Smooth convex optimization with Lipchitz continuous gradient}

Let us assume that a function $f(x)$ is a smooth convex function with $L$-Lipchitz gradient w.r.t. a norm $\norm{\cdot}$, then
\begin{gather}
0 \leq f(x) - f(y) - \langle\nabla f(y), x - y \rangle \leq \frac{L}{2}\norm{x - y}^2  \,\,\,\, \forall x,y \in Q.
\end{gather}
Thus, we have that $\psi_{\delta_k}(x,y) = \langle\nabla f(y), x - y \rangle$ is $(\delta,\allowbreak L,\allowbreak \norm{})$--model of the function $f$, $f_{\delta_k}(y) = f(y)$, and $\delta_k = 0$ for all $k \geq 0$.

\item {\it Convex optimization with Holder continuous subgradients}

Let us assume that a function $f$ is a convex function with Holder continuous subgradients: exists $\nu\in[0,1]$ such that
\begin{equation*}
\norm{\nabla f(x) - \nabla f(y)}_* \leq L_\nu\norm{x - y}^\nu\,\,\,\,\forall x,y \in Q.
\end{equation*}
Then (see \cite{nesterov2015universal})
\begin{gather*}
0 \leq f(x) - f(y) - \langle\nabla f(y), x - y \rangle \leq \frac{L(\delta)}{2}\norm{x - y}^2 + \delta \,\,\,\, \forall x,y \in Q,
\end{gather*}
where \begin{gather*}L(\delta)=L_\nu\left[\frac{L_\nu}{2\delta}\frac{1-\nu}{1+\nu}\right]^\frac{1-\nu}{1+\nu}\end{gather*} and $\delta > 0$ is a controlling parameter. Hence, $\psi_{\delta_k}(x,y) = \langle\nabla f(y), x - y \rangle$ is $(\delta, L(\delta), \norm{})$--model of the function $f$.

\item {\it Composite optimization}

Let us consider the composite optimization problem \cite{nesterov2013gradient}:
\begin{align*}
f(x) := g(x) + h(x) \rightarrow \min_{x \in Q},
\end{align*}
where $g(x)$ is a smooth convex function with $L$--Lipchitz continuous gradients w.r.t. a $\norm{}$ and $h(x)$ is a convex function (not necessarily smooth). For the optimization problem, we have that
\begin{multline*}
0 \leq f(x) - f(y) - \langle\nabla g(y), x - y \rangle - h(x) + h(y) \leq \frac{L}{2}\norm{x - y}^2  \\\quad \forall x,y \in Q.
\end{multline*}
Hence, a function $\psi_{\delta_k}(x,y) = \langle\nabla g(y), x - y \rangle + h(x) - h(y)$ is $(\delta, L(\delta), \norm{})$--model of the function $f$, $f_{\delta_k}(y) = f(y)$, and $\delta_k = 0$ for all $k \geq 0$.
\end{enumerate}

Note that in papers \cite{gasnikov2019fast, stonyakin2019gradient}, we give more examples, including the conditional gradient method (Frank--Wolfe) \cite{ben-tal2015lectures}, superposition of functions \cite{nemirovski1995information, lan2015bundle}, min--min problem \cite{gasnikoveffectivnie}, saddle point problem \cite{gasnikoveffectivnie}. There is an example of $(\delta, L, V)$--model \cite{stonyakin2019gradient} for optimization problem which arises in an electoral model for clustering \cite{nesterov2018soft}.

\subsection{Gradient method}
\label{subsection:gradient_descent}
In \cite{gasnikov2019fast}, the following results were obtained using the concept of inexact model of a function. Consider the adaptive gradient method for the optimization problem \eqref{mainTask}. In Algorithm \ref{algorithm:adpative_gradient_descent} we assume that we have a starting point $x_0$, a local approximation $L_0$ of Lipchitz parameter of a function gradient at a point $x_0$. Also, as the input of the algorithm, we feed sequences $\{\delta_k\}_{k\geq0}$ and $\{\widetilde{\delta}_k\}_{k\geq0}$. 
We assume that on every step $k$, the method has access to $(\delta_k, \bar{L}_{k+1}, \norm{})$--model. In general, a constant $\bar{L}_{k+1}$ can vary from iteration to iteration; we only consider that $(\delta_k, \bar{L}_{k+1}, \norm{})$--model exists. We do not use the constant $\bar{L}_{k+1}$ in Algorithm \ref{algorithm:adpative_gradient_descent} explicitly; furthermore
, our method adapts to this constant. The sequence $\{\widetilde{\delta}_k\}_{k\geq0}$ represents inexact solutions from Definition \ref{solNemirovskiy}, which may be zero, constant, or vary from iteration to iteration in different problems.

\begin{algorithm}
\caption{Adaptive gradient method with $(\delta, L, \norm{})$--model}
\begin{algorithmic}[1]
\label{algorithm:adpative_gradient_descent}
\STATE \textbf{Input:} Starting point $x_0$, sequences $\{\delta_k\}_{k\geq0}$, $\{\widetilde{\delta}_k\}_{k\geq0}$ and $L_0 > 0$.
\STATE $L_1 := \frac{L_0}{2}.$
\FOR{$k \geq 0$}
\STATE Find a minimal integer $i_k\geq 0$ such that
\begin{equation}
\label{algorithm:adpative_gradient_descent:exit}
f_{\delta_k}(x_{k+1}) \leq f_{\delta_k}(x_{k}) + \psi_{\delta_k}(x_{k+1}, x_{k}) +\frac{L_{k+1}}{2}\norm{x_{k} - x_{k+1}}_2^2 + \delta_k,
\end{equation}
where $L_{k+1} = 2^{i_k-1}L_k$, $A_{k+1} := A_k + \frac{1}{L_{k+1}}$.
\begin{equation*}
\phi_{k+1}(x) := \frac{1}{L_{k+1}}\psi_{\delta_k}(x, x_k) + V[x_k](x), \quad
x_{k+1} := {\arg\min_{x \in Q}}^{\widetilde{\delta}_k} \phi_{k+1}(x).
\end{equation*}
\ENDFOR
\end{algorithmic}
\end{algorithm}

In \cite{gasnikov2019fast}, the following convergence rate is derived for Algorithm \ref{algorithm:adpative_gradient_descent}. 

\begin{theorem}[\cite{gasnikov2019fast}]
  \label{theorem:adpative_gradient_descent}
  Let $V[x_0](x_*) \leq R^2$, where $x_0$ is a starting point, $x_*$ is the closest point to $x_0$ in terms of Bregman divergence, a function $f$ is a convex function, for $\delta_k$ and $x_k$ from Algorithm \ref{algorithm:adpative_gradient_descent} we can always find a constant $\bar{L}_{k+1} > 0$ such that $(\delta_k, \bar{L}_{k+1}, \norm{})$--model $\psi_{\delta_k}(\cdot, x_k)$ exists at a point $x_k$, and 
  $\bar{x}_N= \frac{1}{A_N}\sum_{k=0}^{N-1}\alpha_{k+1} x_{k+1}.$
  For Algorithm \ref{algorithm:adpative_gradient_descent}, the following convergence rate holds:
  \begin{align*}
  f(\bar{x}_N) - f(x_*) \leq \frac{R^2}{A_N} + \frac{1}{A_N}\sum_{k=0}^{N-1}\widetilde{\delta}_k + \frac{2}{A_N}\sum_{k=0}^{N-1}\alpha_{k+1}\delta_k.
  \end{align*}
  If we additionally assume that on every step $k$, inexact model $(\delta_k, L, \norm{})$--model exists with a fixed parameter $L$ (in other words, $\bar{L}_{k} \leq L$ for all $k \geq 0$), then
  \begin{align}
  \label{theorem:adpative_gradient_descent:convergence_rate_L}
  f(\bar{x}_N) - f(x_*) \leq \frac{2LR^2}{N} + \frac{2L}{N}\sum_{k=0}^{N-1}\widetilde{\delta}_k + \frac{2}{A_N}\sum_{k=0}^{N-1}\alpha_{k+1}\delta_k.
  \end{align}
\end{theorem}

There are three terms in \eqref{theorem:adpative_gradient_descent:convergence_rate_L} from Theorem \ref{theorem:adpative_gradient_descent}: convergence rate, accumulated error from auxiliary problems, and accumulated error from the inexact model of a function. For the simplicity of the analysis, let us suppose that $\widetilde{\delta}_k = \widetilde{\delta}$ and $\delta_k = \delta$ for all $k \geq 0$, then from \eqref{theorem:adpative_gradient_descent:convergence_rate_L}, we can get a more convenient convergence rate estimate:
  \begin{align}
  \label{theorem:adpative_gradient_descent:convergence_rate_L:simple}
  f(\bar{x}_N) - f(x_*) \leq \frac{2LR^2}{N} + 2L \widetilde{\delta} + \delta.
  \end{align}
From the last inequality, we can conclude that the derived convergence rate corresponds to the convergence rate of the nonaccelerated gradient method \cite{nesterov2004}, while errors $\widetilde{\delta}$ and $\delta$ do not accumulate with the number of algorithm iterations. In Section \ref{subsection:fast_gradient_descent}, we consider the accelerated version of the proposed algorithm, which has a different nature of the convergence rate with respect to $\widetilde{\delta}$ and $\delta$.

Note that we use brute--force search from $0$ to infinity in order to find an integer $i_k$ in Algorithm \ref{algorithm:adpative_gradient_descent}. However, the assumption about the existence of $(\delta_k, L_{k+1}, \norm{})$--model at a point $x_k$ ensures that this process is finite. Moreover, we can show that in ``average'' an integer $i_k$ for which \eqref{algorithm:adpative_gradient_descent:exit} is satisfied is equal $1$ (see \cite{nesterov2015universal}, p. 7--8). Therefore, in ``average'', $\psi_{\delta_k}(\cdot, x_k)$ is requested 2 times in every iteration of Algorithm \ref{algorithm:adpative_gradient_descent}.

\subsection{Fast gradient method}
\label{subsection:fast_gradient_descent}

Let us consider the accelerated version of the algorithm from Section  \ref{subsection:gradient_descent}. In \cite{gasnikov2019fast}, we propose Algorithm \ref{algorithm:fast_adpative_gradient_descent} and prove the corresponding Theorem \ref{theorem:fast_adpative_gradient_descent}.

\begin{algorithm}
\caption{Adaptive fast gradient method with $(\delta, L, \norm{})$--model}
\begin{algorithmic}[1]
\label{algorithm:fast_adpative_gradient_descent}
\STATE \textbf{Input:} Starting point $x_0$, sequences $\{\delta_k\}_{k\geq0}$, $\{\widetilde{\delta}_k\}_{k\geq0}$ and $L_0 > 0$.
\STATE $y_0 := x_0,\,
u_0 := x_0,\,
L_1 := \frac{L_0}{2},\,
\alpha_0 := 0,\,
A_0 := \alpha_0.$
\FOR{$k \geq 0$}
\STATE Find a minimal integer $i_k\geq 0$ such that
\begin{equation*}
f_{\delta_k}(x_{k+1}) \leq f_{\delta_k}(y_{k+1}) + \psi_{\delta_k}(x_{k+1}, y_{k+1})  + \frac{L_{k+1}}{2}\norm{x_{k+1} - y_{k+1}}^2 + \delta_k,
\end{equation*}
where $L_{k+1} = 2^{i_k-1}L_k$, $\alpha_{k+1}$ it the largest root of $A_k + \alpha_{k+1} = L_{k+1}\alpha^2_{k+1}$, $A_{k+1} := A_k + \alpha_{k+1}$.
\begin{gather*}
y_{k+1} := \frac{\alpha_{k+1}u_k + A_k x_k}{A_{k+1}},
\end{gather*}
\begin{equation*}
\begin{gathered}
\phi_{k+1}(x) = V[u_k](x) + \alpha_{k+1}\psi_{\delta_k}(x, y_{k+1}),\\
u_{k+1} := {\arg\min_{x \in Q}}^{\widetilde{\delta}_k}\phi_{k+1}(x),
\end{gathered}
\end{equation*}
\begin{gather*}
x_{k+1} := \frac{\alpha_{k+1}u_{k+1} + A_k x_k}{A_{k+1}}.
\end{gather*}
\ENDFOR
\end{algorithmic}
\end{algorithm}

\begin{theorem}[\cite{gasnikov2019fast}]
  \label{theorem:fast_adpative_gradient_descent}
  Let $V[x_0](x_*) \leq R^2$, where $x_0$ is a starting point, $x_*$ is the closest point to $x_0$ in terms of Bregman divergence, a function $f$ is a convex function and for $\delta_k$ and $y_{k+1}$ from Algorithm \ref{algorithm:fast_adpative_gradient_descent} we can always find a constant $\bar{L}_{k+1} > 0$ such that $(\delta_k, \bar{L}_{k+1}, \norm{})$--model $\psi_{\delta_k}(\cdot, y_{k+1})$ exists at a point $y_{k+1}$.
  For Algorithm \ref{algorithm:fast_adpative_gradient_descent}, the following convergence rate holds:
  \begin{align*}
  f(x_N) - f(x_*) \leq \frac{R^2}{A_{N}} + \frac{\sum_{k=0}^{N-1}\widetilde{\delta}_k}{A_{N}} + \frac{2\sum_{k = 0}^{N-1}\delta_kA_{k+1}}{A_{N}}.
  \end{align*}
  If we additionally assume that on every step $k$, inexact model $(\delta_k, L, \norm{})$--model exists with a fixed parameter $L$ (in other words, $\bar{L}_{k} \leq L$ for all $k \geq 0$), then
  \begin{align}
  \label{theorem:fast_adpative_gradient_descent:convergence_rate_L}
  f(x_N) - f(x_*) \leq \frac{8LR^2}{(N+1)^2} + \frac{8L\sum_{k=0}^{N-1}\widetilde{\delta}_k}{(N+1)^2} + \frac{2\sum_{k = 0}^{N-1}\delta_kA_{k+1}}{A_{N}}.
  \end{align}
\end{theorem}

As in Section \ref{subsection:gradient_descent}, we assume that 
$\widetilde{\delta}_k = \widetilde{\delta}$ and $\delta_k = \delta$ for all $k \geq 0$, then from the inequality \eqref{theorem:fast_adpative_gradient_descent:convergence_rate_L}, we have:
\begin{align*}
  f(x_N) - f(x_*) \leq \frac{8LR^2}{(N+1)^2} + \frac{8L\widetilde{\delta}}{N+1} + N \delta.
\end{align*}
Comparing the last inequality with \eqref{theorem:adpative_gradient_descent:convergence_rate_L:simple} we can conclude that Algorithm \ref{theorem:fast_adpative_gradient_descent} has the convergence rate of the fast gradient method, while an error $\delta$ linearly accumulates with the number of the algorithm iterations. In particular, the impact of an error $\widetilde{\delta}$ decreases. If we compare methods w.r.t. $\widetilde{\delta}$, then Algorithm \ref{algorithm:fast_adpative_gradient_descent} is more effective than Algorithm \ref{algorithm:adpative_gradient_descent}. While the conclusion w.r.t. an error $\delta$ is not so unequivocal and depends on an error $\delta$. More details reader can find in the paper \cite{devolder2011}.

Similarly, as in Section \ref{subsection:gradient_descent}, we can conclude, that in ``average'' an integer $i_k$ is equal to $1$ \cite{nesterov2015universal}.

\subsection{Gradient method with relative smoothness}

Let us consider a simplified version of the algorithm from Section \ref{subsection:gradient_descent}. The following method works with functions supported by an oracle from Definition \ref{definition:model_relative} with relative smoothness. The optimization method from Section \ref{subsection:gradient_descent} (Algorithm \ref{algorithm:adpative_gradient_descent}) is not applicable to numerous optimization problems (see \cite{lu2018}). Further, we consider Algorithm  \ref{algorithm:gradient_descent_relative} and corresponding Theorem \ref{theorem:gradient_descent_relative}.

In this section, we relax the assumption about a prox--function $d$ and replace 1--strong convexity condition with the only convexity of a function $d$. This allows us to apply Theorem \ref{theorem:gradient_descent_relative} in more general cases.

\begin{algorithm}
\caption{Gradient method with $(\delta, L, V)$--model}
\label{algorithm:gradient_descent_relative}
\begin{algorithmic}[1]
\STATE \textbf{Input:} Starting point $x_0$, $L > 0$ and 
$\delta,\widetilde{\delta}>0$.
\FOR{$k \geq 0$}
\STATE \begin{equation}
\phi_{k+1}(x) := \psi_{\delta}(x, x_k)+L V[x_k](x), \quad
x_{k+1} := {\arg\min_{x \in Q}}^{\widetilde{\delta}} \phi_{k+1}(x). 
\end{equation}
\ENDFOR
\end{algorithmic}
\end{algorithm}

\begin{theorem}[\cite{stonyakin2019gradient}]
  \label{theorem:gradient_descent_relative}
  Let $V[x_0](x_*) \leq R^2$, where $x_0$ is a starting point, is the closest point to $x_0$ in terms of Bregman divergence, a function $f$ is a convex function, $(\delta, L, V)$--model $\psi_{\delta}(\cdot, x_k)$ exists for a function $f$ on a set $Q$, and $\bar{x}_N= \frac{1}{N}\sum_{k=0}^{N-1} x_{k+1}.$ 
  For Algorithm \ref{algorithm:gradient_descent_relative}, the following convergence rate holds:
  \begin{align*}
  f(\bar{x}_N) - f(x_*) \leq \frac{LR^2}{N} + \widetilde{\delta} + \delta.
  \end{align*}
\end{theorem}

It would be natural to develop the fast gradient descent with relative smoothness by analogy with Section \ref{subsection:fast_gradient_descent}. However, in general, for optimization problems supported by relative smoothness, convergences rate of nonaccelerated methods can not be improved up to a constant factor (see \cite{dragomir2019optimal}).

\subsection{Primal--dual adaptive gradient method}

\label{subsection:primal_dual_gradient_descent}

In this section, we consider the primal--dual gradient method. The main goal of primal--dual type methods is to find $\varepsilon$--solution of both the primal problem \eqref{mainTask} and the corresponding dual problem. Let us introduce additional assumption on a set $Q$. Consider the following setup for a set $Q$:
  \begin{align}
  \label{primal_dual_Q}
  Q = \{x~|~x \in \widetilde{Q},~f_i(x) \leq 0~ \forall i \in [1, m]\},
  \end{align}
where for all $i$ a function $f_i(x):\widetilde{Q}\rightarrow\mathbb{R}$ is a convex function, and a set $\widetilde{Q}$ is a convex and closed set. Let us define a vector-valued function $F$:
  \begin{align*}
  F(x) = [f_1(x), \dots, f_m(x)]^T.
  \end{align*}
Thus, we rewrite \eqref{mainTask} as
  \begin{align}
  \label{mainTaskPrimal}
  f(x)  \rightarrow \min_{x \in \widetilde{Q},~F(x) \leq 0}.
  \end{align}
Let us construct the Lagrange dual problem. Using a definition
  \begin{align}
  \label{dual_function}
  g(z) = \max_{x \in \widetilde{Q}} [- f(x)-\langle z, F(x)  \rangle].
  \end{align}
we obtain the dual problem for the primal problem \eqref{mainTaskPrimal}:
  \begin{align}
  \label{dual_task_dop}
  g(z)  \rightarrow \min_{z \in \mathbb{R}^{m}_+}.
  \end{align}
From now on, consider the strong duality assumption \cite{boyd2004convex} (p. 226).

The feature to restore a solution of a dual problem is proven to be very useful in various optimization tasks for which it is faster to find an optimal point in a primal problem than in a dual problem. For instance, this property is used in transportation tasks \cite{baymurzina2019universal, gasnikov2018dualtransport, gasnikoveffectivnie}.

\begin{definition}
Let a point $x_*$ be a solution of a primal optimization problem
\begin{equation}
    \label{tmp_primal}
    p(x) \rightarrow \min_{x \in \widetilde{Q},~G(x) \leq 0}.
\end{equation}
A point $z_*$ is a solution of a dual optimization problem
\begin{align*}
  h(z)  \rightarrow \min_{z \in \mathbb{R}^{m}_+},
\end{align*}
for \eqref{tmp_primal}, where $z$ are dual variable w.r.t. constraints $G(x) \leq 0$. We define operator \textit{$\argdualmin$} that depends on functions $p$ and $G$ and returns points $x_*$ and $z_*$:
  \begin{align*}
  (x_*, z_*) := \argdualmin_{x \in \widetilde{Q}}(p(x), G(x)).
  \end{align*}
\end{definition}

In \cite{tyurin2020primal}, we propose Algorithm \ref{algorithm:primal_dual_gradient_descent} and corresponding Theorem \ref{theorem:primal_dual_gradient_descent}.

\begin{algorithm}
\caption{{Primal--dual adaptive gradient method with $(\delta, L, \norm{})$--model}}
\begin{algorithmic}[1]
\label{algorithm:primal_dual_gradient_descent}
\STATE \textbf{Input:} Starting point $x_0$, $L_0 > 0$, and sequence
$\{\delta_k\}_{k\ge0}$.
\STATE $A_0 := 0$
\FOR{$k \geq 0$}
\STATE \label{loop_state} Find a minimal integer $i_k\geq 0$ such that
\begin{equation*}
f_{\delta_k}(x_{k+1}) \leq f_{\delta_k}(x_{k}) + \psi_{\delta_k}(x_{k+1}, x_{k}) + \frac{L_{k+1}}{2}\norm{x_{k+1} - x_{k}}^2 + \delta_k,
\end{equation*}
where $L_{k+1} := 2^{i_k-1}L_k$, $A_{k+1} := A_k + \frac{1}{L_{k+1}}$.
\begin{equation}
\begin{gathered}
\label{primaldualauxilary}
\phi_{k+1}(x) := \psi_{\delta_k}(x, x_k)+L_{k+1}V[x_k](x), \\
(x_{k+1}, z_{k+1}) := \argdualmin_{x \in \widetilde{Q}}(\phi_{k+1}(x), F(x)).
\end{gathered}
\end{equation}
\ENDFOR
\end{algorithmic}
\end{algorithm}

\begin{theorem}[\cite{tyurin2020primal}]
\label{theorem:primal_dual_gradient_descent}
  Let $\bar{x}_N= \frac{1}{A_N}\sum_{k=0}^{N-1}\frac{x_{k+1}}{L_{k+1}}$, $\bar{z}_N= \frac{1}{A_N}\sum_{k=0}^{N-1}\frac{z_{k+1}}{L_{k+1}}$, $V[x_0](\allowbreak x(\bar{z}_N)) \allowbreak \leq R^2$, $x_0$ is a starting point, $x(\bar{z}_N)$ is the maximum point in \eqref{dual_function} with $z = \bar{z}_N$, a function $f$ is a convex function, and for $\delta_k$ and $x_k$ from Algorithm \ref{algorithm:primal_dual_gradient_descent} we can always find a constant $\bar{L}_{k+1} > 0$ such that $(\delta_k, \bar{L}_{k+1}, \norm{})$--model $\psi_{\delta_k}(\cdot, x_k)$ exists at a point $x_k$. For Algorithm \ref{algorithm:primal_dual_gradient_descent}, the following convergence rate holds:
  \begin{align*}
f(\bar{x}_N) + g(\bar{z}_N) \leq \frac{R^2}{A_N} + \frac{1}{A_N}\sum_{k=0}^{N-1}\frac{2\delta_k}{L_{k+1}}.
\end{align*}
\end{theorem}

The theorem fully agrees with results from Theorem  \ref{theorem:adpative_gradient_descent}, taking into account new assumptions about a set $Q$ and condition, that $\widetilde{\delta}_k = 0$ for all $k \geq 0$. However, in Theorem \ref{theorem:primal_dual_gradient_descent}, the convergence rate is proved for a duality gap $f(\bar{x}_N) + g(\bar{z}_N)$.

There are different approaches to restore a dual $\varepsilon$--solution while $\varepsilon$--solution of a primal task is calculating. In a series of our papers \cite{anikin2017dual, dvurechensky2020stable}, dual variables are recovering with the Lagrange function of optimization problem \eqref{dual_task_dop} as a method works. With this approach, conditions in the optimization problem are violated. In an alternative approach \cite{nesterov2009primal}, dual variables are recovering using an auxiliary problem (see, for example, \eqref{primaldualauxilary}); however, worser duality gap bounds can be obtained. Indeed, we have $V[x_0](\allowbreak x(\bar{z}_N))$ instead of $V[x_0](x_*)$. Methods from this and the next section inherit the idea from \cite{nesterov2009primal}.

\subsection{Primal--dual adaptive fast gradient method}

In this section, we consider the accelerated version of the algorithm from Section \ref{subsection:primal_dual_gradient_descent}. Let us study the same assumption on a set $Q$ as in Section \ref{subsection:primal_dual_gradient_descent}. In \cite{tyurin2020primal}, we propose Algorithm \ref{algorithm:primal_dual_fast_gradient_descent} and prove the corresponding Theorem \ref{theorem:primal_dual_fast_gradient_descent}.

% \at{Отметим, что полученный прямо--двойственный быстрый градиентный метод наследует идеи и подходы из работ \cite{anikin2017dual, dvurechensky2020stable}, но при этом обобщается на модельную общность. Также, работа \cite{anikin2017dual} содержит в себе некоторые практические и теоретические замечания, которые применимы и для оптимизационного метода из данного раздела. В работе \cite{dvurechensky2020stable} рассматривается ускоренный прямо--двойственный метод для решения транспортных задач и показывается, что полученный метод на транспортных задачах с маленьким параметром регуляризации работает на практике быстрее, чем метод Синхорна \cite{sinkhorn1974diagonal}. Важным замечанием является то, что алгоритм 2 из \cite{dvurechensky2020stable} является частным случаем алгоритма \ref{algorithm:primal_dual_fast_gradient_descent}, поэтому некоторые полученные выводы из \cite{dvurechensky2020stable} можно переносить на алгоритм \ref{algorithm:primal_dual_fast_gradient_descent}, включая его практическую значимость для транспортных задач.}

\begin{algorithm}
\caption{Primal--dual adaptive fast gradient method with $(\delta, L, \norm{})$--model}
\begin{algorithmic}[1]
\label{algorithm:primal_dual_fast_gradient_descent}
\STATE \textbf{Input:} Starting point $x_0$, sequence $\{\delta_k\}_{k\geq0}$ and $L_0 > 0$.
\STATE $y_0 := x_0,\,
u_0 := x_0,\,
L_1 := \frac{L_0}{2},\,
\alpha_0 := 0,\,
A_0 := \alpha_0.$
\FOR{$k \geq 0$}
\STATE Find a minimal integer $i_k\geq 0$ such that
\begin{equation*}
f_{\delta_k}(x_{k+1}) \leq f_{\delta_k}(y_{k+1}) + \psi_{\delta_k}(x_{k+1}, y_{k+1})  + \frac{L_{k+1}}{2}\norm{x_{k+1} - y_{k+1}}^2 + \delta_k,
\end{equation*}
where $L_{k+1} = 2^{i_k-1}L_k$, $\alpha_{k+1}$ is the largest root of $A_k + \alpha_{k+1} = L_{k+1}\alpha^2_{k+1}$, $A_{k+1} := A_k + \alpha_{k+1}$.
\begin{gather*}
y_{k+1} := \frac{\alpha_{k+1}u_k + A_k x_k}{A_{k+1}},
\end{gather*}
\begin{equation*}
\begin{gathered}
\phi_{k+1}(x) := \psi_{\delta_k}(x, y_{k+1})+L_{k+1}V[u_k](x), \\
(x_{k+1}, z_{k+1}) := \argdualmin_{x \in \widetilde{Q}}(\phi_{k+1}(x), F(x)).
\end{gathered}
\end{equation*}
\begin{gather*}
x_{k+1} := \frac{\alpha_{k+1}u_{k+1} + A_k x_k}{A_{k+1}}.
\end{gather*}
\ENDFOR
\end{algorithmic}
\end{algorithm}

\begin{theorem}[\cite{tyurin2020primal}]
\label{theorem:primal_dual_fast_gradient_descent}
  Let $\bar{z}_N= \frac{1}{A_N}\sum_{k=0}^{N-1}\alpha_{k+1}z_{k+1}$, $V[x_0](x(\bar{z}_N)) \leq R^2$, $x_0$ is a starting point, $x(\bar{z}_N)$ is the maximum point in \eqref{dual_function} with $z = \bar{z}_N$, a function $f$ is a convex function, and for $\delta_k$ and $y_{k+1}$ from Algorithm \ref{algorithm:primal_dual_fast_gradient_descent} we can always find a constant $\bar{L}_{k+1} > 0$ such that $(\delta_k, \bar{L}_{k+1}, \norm{})$--model $\psi_{\delta_k}(\cdot, y_{k+1})$ exists at a point $y_{k+1}$. For Algorithm \ref{algorithm:primal_dual_fast_gradient_descent}, the following convergence rate holds:
%   \begin{align*}
%   f(x_{N})\leq \min_{x \in \widetilde{Q}}\left[\frac{1}{A_N}\sum_{k=0}^{N-1}\alpha_{k+1}\Big(f_{\delta_k}(y_{k+1}) + \psi_{\delta_k}(x,y_{k+1})\Big) + \langle \bar{z}_N, F(x) \rangle + \frac{V(x, u_0)}{A_N}\right] + \frac{2}{A_N}\sum_{k=0}^{N-1}A_{k+1}\delta_k
% \end{align*}
% и 
  \begin{align*}
f(x_{N}) + g(\bar{z}_N) &\leq \frac{R^2}{A_N} + \frac{2}{A_N}\sum_{k=0}^{N-1}A_{k+1}\delta_k.
\end{align*}
\end{theorem}

With the assumption \eqref{primal_dual_Q} about a set $Q$, the derived convergence rate fully agrees with the convergence rate from Theorem \ref{theorem:fast_adpative_gradient_descent} if we consider, that $\widetilde{\delta}_k = 0$ for all $k \geq 0$.

\subsection{Stochastic fast gradient method}
\label{subsection:stochastic_gradient_descent}

Let us consider $(\delta_1, \delta_2, L, \norm{})$--model from Definition \ref{definition:model_stochastic} that generalizes $(\delta, L,\allowbreak \norm{}\allowbreak)$--model of a function. Similarly to Section \ref{subsection:gradient_descent} and \ref{subsection:fast_gradient_descent}, in paper \cite{dvinskikh2020accelerated}, we provide convergence rates for methods that work with $(\delta_1, \delta_2, L, \norm{})$--model. One of the most important consequences is that this concept is surprisingly well--suited for stochastic optimization problems \cite{lanlectures, devolder2013exactness}. 
In Algorithm \ref{algorithm:stochastic_fast_adpative_gradient_descent}, we provide the fast gradient method with $(\delta_1, \delta_2, L, \norm{}_2)$--model.

\begin{algorithm}
\caption{Fast gradient method with $(\delta_1, \delta_2, L, \norm{}_2)$--model}
\begin{algorithmic}[1]
\label{algorithm:stochastic_fast_adpative_gradient_descent}
\STATE \textbf{Input:} Starting point $x_0$ and $L > 0$.
\STATE $y_0 := x_0,\,
u_0 := x_0,\,
\alpha_0 := 0,\,
A_0 := \alpha_0.$
\FOR{$k \geq 0$}
\STATE Constant $\alpha_{k+1}$ is the largest root of $A_k + \alpha_{k+1} = L \alpha^2_{k+1}$, $A_{k+1} := A_k + \alpha_{k+1}$.
\begin{gather*}
y_{k+1} := \frac{\alpha_{k+1}u_k + A_k x_k}{A_{k+1}},
\end{gather*}
\begin{equation*}
\begin{gathered}
\phi_{k+1}(x) = \frac{1}{2}\norm{x - u_k}_2^2 + \alpha_{k+1}\psi_{\delta_k}(x, y_{k+1}),\\
u_{k+1} := \arg\min_{x \in Q} \phi_{k+1}(x),
\end{gathered}
\end{equation*}
\begin{gather*}
x_{k+1} := \frac{\alpha_{k+1}u_{k+1} + A_k x_k}{A_{k+1}}.
\end{gather*}
\ENDFOR
\end{algorithmic}
\end{algorithm}

Now, we formulate the convergence rate theorem for $(\delta_1, \delta_2, L, \norm{})$--model.

\begin{theorem}[\cite{dvinskikh2020accelerated}]
\label{theorem:stochatic_fast_gradient_descent}
  Let $\frac{1}{2}\norm{x_* - x_0}_2^2 \leq R^2$, where $x_0$ is a starting point, $x_*$ is the closest point to $x_0$ in terms of euclidean distance, a function $f$ is a convex function, $(\delta_1^k, \delta_2^k, L, \norm{}_2)$--model $\psi_{\delta_k}(\cdot, y_{k+1})$ exists at a point $y_{k+1}$ from Algorithm \ref{algorithm:stochastic_fast_adpative_gradient_descent}.
  For Algorithm \ref{algorithm:stochastic_fast_adpative_gradient_descent}, the following convergence rate holds:
    \begin{align}
\begin{split}
    \label{theorem:stochatic_fast_gradient_descent:convergence_rate_L}
    f(x_{N}) - f(x_*) &\leq \frac{4LR^2}{N^2}  + \frac{1}{A_N}\sum_{k=0}^{N-1}A_{k}\delta_1^k(x_k,y_{k+1}) \\
    &\hspace{2em}+ \frac{1}{A_N}\sum_{k=0}^{N-1}\alpha_{k+1}\delta_1^k(x_*,y_{k+1}) + \frac{1}{A_N}\sum_{k=0}^{N-1}A_{k+1}\delta_2^k.
    \end{split}
    \end{align}
\end{theorem}

If we additionally assume that  
$\{\delta_1^k\}_{k=0}^{N-1}$ and $\{\delta_2^k\}_{k=0}^{N-1}$ are random sequences with assumptions:
\begin{assumption}
\label{assumption:delta}
Given two sequences $\delta_1^k(y,x)$ and $\delta_2^k$ ($k \geq 0$). Assume that
\begin{center}
$\E \left[\delta_1^k(y,x) |\delta_{1,2}^{k-1},\delta_{1,2}^{k-2},... \right]= 0$,  (conditionally unbiased)
\end{center}
 $\delta_1^k(y,x)$ has  $\left(\hat{\delta}_1\right)^2$--subgaussion conditional variance, $\sqrt{\delta_2^k}$ has $\hat{\delta}_2$--subgaussion conditional second moment.
\end{assumption}

\begin{assumption}
\label{assumption:delta_4}
Given two sequences $\delta_1^k(x,y)$ and $\delta_2^k$ ($k \geq 0$). The random variable $\delta_1^k(x, y)$ has  $\left(\hat{\delta}_1^k(x - y)\right)^2$--subgaussion conditional moment ($\hat{\delta}_1^k(\cdot)$ is a non--random function of one variable) such that 
\begin{enumerate}
    % \item $\hat{\delta}_1^k(x, y) = \hat{\delta}_1^k(x - y)$ for all $x, y \in Q$.
    \item $\hat{\delta}_1^k(\alpha z) \leq \alpha \hat{\delta}_1^k(z)$ for all $\alpha \geq 0$ and $z \in B(0, R)$.
    \item $\hat{\delta}_1 < +\infty$, where $\hat{\delta}_1 \geq \sup_{z \in B(0, R)} \hat{\delta}_1^k(z)$.
\end{enumerate}
\end{assumption}

With high probability \footnote{it means that with probability $\ge 1 - \gamma$, and $\tilde{O}(\cdot)$ means the same as $O(\cdot)$; however, a constant factor depends on $\ln \left( 1/\gamma\right)$.}
\begin{equation*}
f(x_N) - f(x_*) = \tilde{O}\left(\frac{LR^2}{N^2} + \frac{\hat{\delta}_1}{\sqrt{N}} + N \hat{\delta}_2\right).
\end{equation*}
Moreover,
\begin{equation*}
\mathbb{E}[f(x_N)] - f(x_*) = O\left(\frac{LR^2}{N^2} +  N \hat{\delta}_2\right).
\end{equation*}

The stochastic optimization problem is an important case that can be described by $(\delta_1, \delta_2, L, \norm{})$--model. Let us consider the following optimization task:
\begin{equation}
\label{stochastic_problem}
f(x) = \mathbb{E}[f(x,\xi)] \rightarrow \min_{x\in Q},
\end{equation} 
where a set $Q$ is a convex and closed set, $\xi$ is a random variable, the expected value $\mathbb{E}[f(x,\xi)]$ is well--defined and finite for all $x \in Q$,  a function $f$ is a convex function with $L$--Lipschitz continuous gradient, $\nabla f(y,\xi)$ has subgaussian distribution with subgaussian variance $\sigma^2$. For optimization tasks \eqref{stochastic_problem}, we can take $\psi_{\delta}(x,y) = \langle \nabla f(y,\xi), x-y \rangle$, and we can show (see. \cite{dvinskikh2020accelerated}) that for sequences $\{\delta_k^1\}_{k=0}^{N-1}$ and $\{\delta_k^2\}_{k=0}^{N-1}$, the following bounds hold: $\hat{\delta}_1 = O(\sigma R)$ and $\hat{\delta}_2 = O(\sigma^2/L)$. The optimal convergence rate can be obtained for the task \eqref{stochastic_problem} with the help of a mini--batch technique (see \cite{dvinskikh2020accelerated}).

Note that the same reasoning can be applied to composite and min--max optimization tasks.

\subsection{Heuristic adaptive stochastic fast gradient method}
In \cite{ogaltsov2020heuristic}, we propose the heuristic adaptive stochastic fast gradient method based on the adaptive fast gradient method (Algorithm \ref{algorithm:fast_adpative_gradient_descent}) and the nonadaptive stochastic gradient method \cite{dvinskikh2020accelerated} (Algorithm \ref{algorithm:stochastic_fast_adpative_gradient_descent}). For now, it is an open question, if it is possible to add adaptivity to the stochastic fast gradient method in order to preserve convergence rate estimates. Various attempts were made in \cite{bach2019universal, vaswani2019painless, ward2019adagrad, deng2018optimal, levy2018online, iusem2019variance}. A more detailed analysis reader can find in \cite{ogaltsov2020heuristic}. Let us define a mini--batch of gradients as
\begin{gather*}
\widetilde{\nabla}^{m_{k+1}} f(y) = \frac{1}{m_{k+1}}\sum_{j=1}^{m_{k+1}}\nabla f(y;\xi_j),
\end{gather*}
and mini--batch of functions values as
\begin{gather*}
f^{m_{k+1}}(y) = \frac{1}{m_{k+1}}\sum_{j=1}^{m_{k+1}} f(y;\xi_j),
\end{gather*}
where $\xi_j$ are random variables ($j=1,...,m_{k+1}$), $\nabla f(y;\xi_j)$ and $f(y;\xi_j)$ are unbiased estimates of $\nabla f(y)$ and $f(y)$, and $m_{k+1}$ is the number of elements in the mini--batch. In Algorithm \ref{algorithm:fast_stcohastic_adpative_gradient_descent}, we present the heuristic method (see \cite{ogaltsov2020heuristic}) that works with stochastic gradients and stochastic function values. In \cite{dvinskikh2020accelerated}, using an inexact model from Definition  \ref{definition:model_stochastic}, we have the theorem that shows the convergence rate estimate for the nonadaptive version of Algorithm \ref{algorithm:fast_stcohastic_adpative_gradient_descent}.

Note that in paper \cite{ogaltsov2020heuristic}, we have the approbation of our algorithm with the help of experiments. Using practical machine learning tasks MNIST \cite{lecun-mnisthandwrittendigit-2010} and CIFAR \cite{cifar10}, we show that our algorithm convergence faster than popular optimization methods Adam \cite{adam} and AdaGrad \cite{duchi2011adaptive} with logistic regression loss function and linear, neural network, and convolutional neural network backbones.

\begin{algorithm}
\caption{Heuristic adaptive stochastic fast gradient method}
\begin{algorithmic}[1]
\label{algorithm:fast_stcohastic_adpative_gradient_descent}
\STATE \textbf{Input:} Starting point $x_0$, constants $\epsilon > 0$, $L_0 > 0$ and $\sigma_0^2 > 0$.
\STATE $y_0 := x_0,\,
u_0 := x_0,\,
L_1 := \frac{L_0}{2},\,
\alpha_0 := 0,\,
A_0 := \alpha_0.$
\FOR{$k \geq 0$}
\STATE Find a minimal integer $i_k\geq 0$ such that
\begin{equation*}
\begin{gathered}
f^{m_{k+1}}(x_{k+1}) \leq f^{m_{k+1}}(y_{k+1}) + \langle \widetilde{\nabla}^{m_{k+1}} f(y_{k+1}), x_{k+1} - y_{k+1} \rangle +\\  + \frac{L_{k+1}}{2}\norm{x_{k+1} - y_{k+1}}^2 + \frac{\epsilon}{L_{k+1}\alpha_{k+1}},
\end{gathered}
\end{equation*}
where $L_{k+1} = 2^{i_k-1}L_k$, $\widetilde{\alpha}_{k+1}$ it the largest root of $A_k + \alpha_{k+1} = L_{k}\alpha^2_{k+1}$, $\alpha_{k+1}$ it the largest root of $A_k + \alpha_{k+1} = L_{k+1}\alpha^2_{k+1}$, $A_{k+1} := A_k + \alpha_{k+1}$, $m_{k+1} := \left\lceil\frac{3\sigma_0^2\tilde{\alpha}_{k+1}}{\epsilon}\right\rceil$. If $i_k = 0$, then generate i.i.d. $\xi_j$ ($j=1,...,m_{k+1}$).
\begin{gather*}
y_{k+1} := \frac{\alpha_{k+1}u_k + A_k x_k}{A_{k+1}},
\end{gather*}
\begin{equation*}
\begin{gathered}
\hspace{-1.0cm}\phi_{k+1}(x) := \frac{1}{2}\norm{x - u_k}^2_2 + \alpha_{k+1}\left(f^{m_{k+1}}(y_{k+1}) + \langle \widetilde{\nabla}^{m_{k+1}} f(y_{k+1}), x - y_{k+1} \rangle\right),\\
u_{k+1} := \arg\min_{x \in Q}\phi_{k+1}(x),
\end{gathered}
\end{equation*}
\begin{gather*}
x_{k+1} := \frac{\alpha_{k+1}u_{k+1} + A_k x_k}{A_{k+1}}.
\end{gather*}
\ENDFOR
\end{algorithmic}
\end{algorithm}

\section{Conclusion}

This thesis is based on published papers \cite{gasnikov2019fast, stonyakin2019gradient, tyurin2020primal, dvinskikh2020accelerated, ogaltsov2020heuristic, dvurechensky2020stable, anikin2017dual, dvinskikh2020sliding}. 

In papers \cite{gasnikov2019fast, stonyakin2019gradient, tyurin2020primal, ogaltsov2020heuristic}, we developed optimization methods that exploit the concepts of inexact model. Also, we demonstrate various examples of optimization tasks supported by suitable inexact models. In addition to standard tasks from structural optimization like smooth optimization, composite optimization, optimization with Holder continuous subgradients, the proposed concepts of inexact model can describe transportation tasks \cite{gasnikov2019fast, gasnikoveffectivnie}, optimization problem which arises in an electoral model for clustering \cite{stonyakin2019gradient, nesterov2018soft}, etc.

Papers \cite{ogaltsov2020heuristic, dvurechensky2020stable, anikin2017dual, dvinskikh2020sliding} are milestones from the view of the development of the concept of inexact model. Moreover, they motivate further research.

Let us list the main results that are obtained in this thesis and submitted for defense.
\begin{enumerate}
    \item Various concepts of inexact model are developed for gradient-type methods. As shown in the thesis, these concepts can represent a significant number of modern optimization problems.
    \item Different gradient methods are proposed which support concepts of inexact model. For methods from sections \ref{subsection:gradient_descent}--\ref{subsection:stochastic_gradient_descent}, we proved corresponding convergence rate theorems and performed analysis.
   \item The heuristic adaptive stochastic fast gradient method is developed and justified.
    \item Theoretical analysis of primal-dual methods for problems with strongly convex functionals with a simple structure under affine constraints, problems that calculate regularized optimal transport, and problems with the concept of inexact model is carried out.
\end{enumerate}

It is worthwhile to say that some ideas and several examples of optimization problems that are well-described by the concept of inexact model are not listed: 
\begin{enumerate}
    \item
In further research, we are planning to develop methods that work with strongly convex functions \cite{stonyakin2020inexact}.

    \item
In \cite{dvinskikh2020sliding}, we consider the practical optimization task with an objective function that has the form of a sum of smooth strongly convex functions with a smooth regularizer. In this paper, we propose an approach that allows us to derive the optimal bound for the case when the composite part of an objective function is not proximal--friendly.\footnote{A proximal--friendly function is a function that, with a quadratic function, can be ``simply'' minimized.}. In further research, we are planning to generalize this result with the concept of inexact model.

    \item
We are planning to combine the concept of inexact model with coordinate descent methods \cite{dvurechensky1707randomized, nesterov2012efficiency}. As in general stochastic optimization, coordinate descent methods admit efficient optimization steps by calculating a descent direction via a subset of coordinates.

    \item
In the end of Section \ref{subsection:primal_dual_gradient_descent}, we mention the described in \cite{anikin2017dual, dvurechensky2020stable} approach of recovering of dual variables. It would be useful to generalize methods from \cite{anikin2017dual, dvurechensky2020stable} using the concept of inexact model.

\end{enumerate}

% \printbibliography[heading=bibnumbered]
\bibliographystyle{ieeetr}
\bibliography{main.bib}

\begin{thebibliography}{10}

\bibitem{kantarovich1939}
{\em Kantorovich L. \normalfont Mathematical methods of organizing and planning
  pro\-du\-ction // Management Science. 1960. V. 6, №. 4. P. 366--422.}

\bibitem{polyakhistory}
{\em Polyak B. \normalfont History of mathematical programming in the USSR:
  ana\-ly\-zing the phenomenon // Math. Program. 2002. V. 91. № 3. P.
  401--416.}

\bibitem{nemirovkii1979}
{\em Nemirovski A., Yudin D. \normalfont Problem complexity and method
  efficiency in optimization. Wiley--Interscience. 1983.}

\bibitem{nesterov1983}
{\em Nesterov Yu. \normalfont A method of solving a convex programming problem
  with convergence rate $O (1/k^2)$ // Dokl. Akad. Nauk SSSR. 1983. V. 269. №
  3. P. 543--547.}

\bibitem{nesterov2013}
{\em Nesterov Yu. \normalfont Gradient methods for minimizing composite
  functions // Math. Program. 2013. V. 140, № 1. P. 125–161.}

\bibitem{nesterov2015universal}
{\em Nesterov Yu. \normalfont Universal gradient methods for convex
  optimization pro\-blems // Math. Program. 2015. V. 152. № 1--2. P.
  381--404.}

\bibitem{nemirovski1995information}
{\em Nemirovski A. \normalfont Information-based complexity of convex
  programming. Technion. 1995.}

\bibitem{lan2015bundle}
{\em Lan G. \normalfont Bundle-level type methods uniformly optimal for smooth
  and nonsmooth convex optimization // Math. Program. 2015. V. 149. № 1--2.
  P. 1--45.}

\bibitem{baymurzina2019universal}
{\em Baimurzina D., Gasnikov A., Gasnikova E., Kubentaeva M.,
  La\-gu\-nov\-skaya A., Dvurechensky P., Ershov E. \normalfont Universal
  Method of Searching for Equilibria and Stochastic Equilibria in
  Transportation Networks // Computational Mathematics and Mathematical
  Physics. 2019. V. 59. № 1. P. 19--33.}

\bibitem{gasnikov2018dualtransport}
{\em Gasnikov A., Gasnikova E., Nesterov Yu. \normalfont Dual methods for
  finding equilibriums in mixed models of flow distribution in large
  transportation networks // Computational Mathematics and Mathematical
  Phys\-ics. 2017. V. 58. № 9. P. 1395--1403.}

\bibitem{gasnikoveffectivnie}
{\em Gasnikov A. \normalfont Effective numerical methods for finding
  equilibrium in large transport networks. PhD thesis. MFTI, 2016.}

\bibitem{nesterov2018soft}
{\em Nesterov Yu. \normalfont Soft clustering by convex electoral model // CORE
  Dis\-cussion paper. 2018/01. 20p. URL:
  \url{https://alfresco.uclouvain.be/alfresco/service/guest/streamDownload/workspace/SpacesStore/ff42ec88-4339-4223-b05d-b768c71ef4e6/coredp2018_01web.pdf?guest=true}}.

\bibitem{nesterov2004}
{\em Nesterov Yu. \normalfont Lectures on convex optimization. Springer. 2018.}

\bibitem{devolder2011}
{\em Devolder O., Glineur F., Nesterov Yu. \normalfont First-order methods of
  smooth convex optimization with inexact oracle // Math. Program. 2014. V.
  146. № 1--2. P. 37--75.}

\bibitem{nesterov2005}
{\em Nesterov Yu. \normalfont Smooth minimization of non-smooth functions //
  Math. Program. 2005. V. 103. № 1. P. 127--152.}

\bibitem{nesterov2005excessive}
{\em Nesterov Yu. \normalfont Excessive gap technique in nonsmooth convex
  mi\-ni\-mi\-zation // SIAM J. Optimizat. 2005. V. 16, № 1. P. 235--249.}

\bibitem{nesterov2007smoothing}
{\em Nesterov Yu. \normalfont Smoothing technique and its applications in
  semidefinite optimization // Math. Program. 2007. V. 110. № 2. P.
  245--259.}

\bibitem{lemarechal1997}
{\em Lemarechal C., Sagaztizabal C. \normalfont Practice aspects of
  Moreau–Yosida regu\-larization: Theoretical preliminaries // SIAM J.
  Optimizat. 1997. V. 7, № 2. P. 367--385.}

\bibitem{aspremont}
{\em D’Aspremont A, \normalfont Smooth optimization with approximate gradient
  // SIAM J. Optimizat. 2019. V. 19, № 3. P. 1171--1183.}

\bibitem{anikin2017dual}
{\em Anikin A., Gasnikov A., Dvurechensky P., Tyurin A., Chernov A. \normalfont
  Dual approaches to the minimization of strongly convex functionals with a
  simple structure under affine constraints // Computational Mathematics and
  Mathematical Physics. 2017. V. 57. № 8. P. 1262--1276.}

\bibitem{boyd2004convex}
{\em Boyd S., Vandenberghe L. \normalfont Convex optimization. Cambridge
  University Press. 2004.}

\bibitem{nesterov2015complexity}
{\em Nesterov Yu. \normalfont Complexity bounds for primal-dual methods
  minimizing the model of objective function // Mathematical Programming. 2018.
  V. 171, №. 1--2. P. 311--330.}

\bibitem{nesterov2009primal}
{\em Nesterov Yu. \normalfont Primal-dual subgradient methods for convex
  problems. // Mathematical Programming. 2009. V. 120, № 1. P. 221--259.}

\bibitem{goodfellow2016deep}
{\em Goodfellow I., Bengio Y., Courville A. \normalfont Deep learning. MIT
  press. 2016.}

\bibitem{krizhevsky2012imagenet}
{\em Krizhevsky A., Sutskever I., Hinton G. \normalfont Imagenet classification
  with deep convolutional neural networks // In Advances in neural
  in\-for\-ma\-tion processing systems. 2012. P. 1097--1105.}

\bibitem{bregman1967relaxation}
{\em Bregman L. \normalfont The relaxation method of finding the common points
  of convex sets and its application to the solution of problems in convex
  programming // USSR Computational Mathematics and Mathematical Phy\-sics.
  1967. V. 7. № 3. P. 200--217.}

\bibitem{ben-tal2015lectures}
{\em Ben-Tal A., Nemirovski A. \normalfont Lectures on Modern Convex
  Optimization. Philadelphia: SIAM, 2015. URL:
  \url{http://www2.isye.gatech.edu/~nemirovs/Lect_ModConvOpt.pdf}}.

\bibitem{gasnikov2019fast}
{\em Gasnikov A., Tyurin A. \normalfont Fast gradient descent for convex
  minimization problems with an oracle producing a $(\delta, L)$-model of
  function at the requested point // Computational Mathematics and Mathematical
  Phy\-sics. 2019. V. 59. № 7. P. 1085--1097.}

\bibitem{stonyakin2020inexact}
{\em Stonyakin F., Tyurin A., Gasnikov A., Dvurechensky P., Agafonov A.,
  Dvinskikh D., Pasechnyuk D., Artamonov S., Piskunova V. \normalfont Inexact
  relative smoothness and strong convexity for optimization and variat\-ional
  inequalities by inexact model // e-print. arXiv:2001.09013. 2020.}

\bibitem{tyurin2020primal}
{\em Tyurin A. \normalfont Primal--dual fast gradient method with a model //
  Computer Research and Modeling. 2020. V. 12, № 2. P. 263--274.}

\bibitem{bauschke2016}
{\em Bauschke H., Bolte J., Teboulle M. \normalfont A descent lemma beyond
  lipschitz gradient continuity: first-order methods revisited and applications
  // Mathematics of Operations Research. 2016. V. 42. № 2. P. 330--348.}

\bibitem{lu2018}
{\em Lu H., Freund R., Nesterov Yu. \normalfont Relatively smooth convex
  optimization by first-order methods, and applications // SIAM J. Optimizat.
  2018. V. 28, № 1. P. 333--354}.

\bibitem{stonyakin2019gradient}
{\em Stonyakin F., Dvinskikh D., Dvurechensky P., Kroshnin A., Kuznetsova O.,
  Agafonov A., Gasnikov A., Tyurin A., Uribe С., Pasechnyuk D., Artamonov S.
  \normalfont Gradient methods for problems with inexact model of the objective
  // Lecture Notes in Computer Science. 2019. V. 11548. P. 97--114.}

\bibitem{dvinskikh2020accelerated}
{\em Dvinskikh D., Tyurin A., Gasnikov A., Omelchenko S. \normalfont
  Accelerated and nonaccelerated stochastic gradient descent with model
  conception // Mathematical Notes. 2020. V. 108. № 4. In press.}

\bibitem{nesterov2013gradient}
{\em Nesterov Yu. \normalfont Gradient methods for minimizing composite
  functions // Math. Program. 2013. V. 140, № 1. P. 125–161.}

\bibitem{dragomir2019optimal}
{\em Dragomir R., Taylor A., D'Aspremont A., Bolte J. \normalfont Optimal
  complexity and certification of Bregman first-order methods // e-print.
  arXiv: 1911.08510. 2019.}

\bibitem{dvurechensky2020stable}
{\em Dvurechensky P., Gasnikov A., Omelchenko A., Tyurin A. \normalfont A
  stable alternative to Sinkhorn’s algorithm for regularized optimal
  transport // Lecture Notes in Computer Science. 2020. V. 12095. P. 406--423.}

\bibitem{lanlectures}
{\em Lan G. \normalfont Lectures on optimization. Methods for Machine Learning
  // e-print. 2019. URL:
  \url{http://pwp.gatech.edu/guanghui-lan/wp-content/uploads/sites/330/2019/08/LectureOPTML.pdf}}.

\bibitem{devolder2013exactness}
{\em Devolder O. \normalfont Exactness, inexactness and stochasticity in
  first-order meth\-ods for large-scale convex optimization. PhD thesis. CORE
  UCL, 2013.}

\bibitem{ogaltsov2020heuristic}
{\em Ogaltsov A., Tyurin A. \normalfont A heuristic adaptive fast gradient
  method in stochastic optimization problems // Computational Mathematics and
  Math\-ematical Physics. 2019. V. 60. № 7. P. 1108--1115}.

\bibitem{bach2019universal}
{\em Bach F., Levy K.Y. \normalfont A universal algorithm for variational
  inequalities adaptive to smoothness and noise // COLT, 2019.}

\bibitem{vaswani2019painless}
{\em Vaswani S., Mishkin A., Laradji I., Schmidt M., Gidel G., Lacoste-Julien
  S. \normalfont Painless Stochastic Gradient: interpolation, line-search, and
  convergence rates // NIPS, 2019.}

\bibitem{ward2019adagrad}
{\em Ward R., Wu X., Bottou L. \normalfont AdaGrad stepsizes: sharp convergence
  over nonconvex landscapes, from any initialization // ICML, 2019.}

\bibitem{deng2018optimal}
{\em Deng Q., Cheng Y., Lan G. \normalfont Optimal adaptive and accelerated
  stochastic gradient descent // e-print. arXiv:1810.00553. 2018.}

\bibitem{levy2018online}
{\em Levy K.Y., Yurtsever A., Cevher V. \normalfont Online adaptive methods,
  uni\-ver\-sa\-lity and acceleration // NIPS, 2018.}

\bibitem{iusem2019variance}
{\em Iusem A.N., Jofre A., Oliveira R.I., Thompson P. \normalfont
  Variance-based extra\-gradient methods with line search for stochastic
  variational inequalities // SIAM J. Optimizat. 2019. V. 29, № 1. P.
  175--206}.

\bibitem{lecun-mnisthandwrittendigit-2010}
{\em LeCun Y., Bottou L., Bengio Y., Haffner P. \normalfont Gradient-based
  learning applied to document recognition // Proceedings of the IEEE. 1998. V.
  86. № 11. P. 2278--2324.}

\bibitem{cifar10}
{\em Krizhevsky A. \normalfont Learning Multiple Layers of Features from Tiny
  Images. PhD thesis. University of Toronto, 2009.}

\bibitem{adam}
{\em Kingma D.P., Ba J. \normalfont Adam: a method for stochastic optimization
  // ICLR, 2015.}

\bibitem{duchi2011adaptive}
{\em Duchi J., Hazan E., Singer Y. \normalfont Adaptive subgradient methods for
  online learning and stochastic optimization // Journal of Machine Learning
  Research. 2011. V. 12. № Jul. P. 2121--2159.}

\bibitem{dvinskikh2020sliding}
{\em Dvinskikh D., Omelchenko A., Gasnikov A., Tyurin A. \normalfont
  Accelerated gradient sliding for minimizing the sum of functions // Doklady
  Mathe\-matics. 2020. V. 101. № 3. In press.}

\bibitem{dvurechensky1707randomized}
{\em Dvurechensky P., Gasnikov A., Tyurin A. \normalfont Randomized similar
  triangles method: a unifying framework for accelerated randomized
  optimization methods (coordinate descent, directional search,
  derivative--free method) // e-print. arXiv: 1707.08486. 2017.}

\bibitem{nesterov2012efficiency}
{\em Nesterov Yu. \normalfont Efficiency of coordinate descent methods on
  huge--scale optimization problems // SIAM J. Optimizat. 2012. V. 22, № 2.
  P. 341--362.}

\end{thebibliography}

\iftoggle{thesis}{

\newpage

\section{Приложение}

\appendix

\section{Статья ``Fast gradient descent for convex mi\-ni\-mi\-zation problems with an oracle producing a $(\delta, L)$-model of function at the requested point''}

\hspace{1cm}

\textbf{Авторы:} Gasnikov A., Tyurin A.

\textbf{Аннотация:} A new concept of $(\delta, L)$--model of a function that is a generalization of the Devolder–Glineur–Nesterov $(\delta, L)$--oracle is proposed. Within this concept, the gradient descent and fast gradient descent methods are constructed and it is shown that constructs of many known methods (composite methods, level methods, conditional gradient and proxi\-mal meth\-ods) are particular cases of the methods proposed in this paper.

\includepdf[pages=-]{CMMP1085.pdf}

\section{Статья ``Gradient methods for problems with in\-exact model of the objective''}

\hspace{1cm}

\textbf{Авторы:} Stonyakin F., Dvinskikh D., Dvurechensky P., Kroshnin A., Kuznetsova O., Agafonov A., Gasnikov A., Tyurin A., Uribe С., Pasechnyuk D., Artamonov S.

\textbf{Аннотация:} We consider optimization methods for convex minimization problems under inexact information on the objective function. We introduce inexact model of the objective, which as a particular cases includes inexact oracle [16] and relative smoothness condition [36]. We analyze gradient meth\-od which uses this inexact model and obtain convergence rates for convex and strongly convex problems. To show potential applications of our general framework we consider three particular problems. The first one is clustering by electorial model introduced in [41]. The second one is approximating optimal transport distance, for which we propose a Proximal Sinkhorn algo\-rithm. The third one is devoted to approximating optimal transport barycen\-ter and we propose a Proximal Iterative Bregman Projections algorithm. We also illustrate the practical performance of our algorithms by numerical experiments.

\includepdf[pages=-]{Stonyakin2019_Chapterh.pdf}

\section{Статья ``Primal-dual fast gradient method with a model''}

\hspace{1cm}

\textbf{Авторы:} Tyurin A.

\textbf{Аннотация:} In this work we consider a possibility to use the conception of $(\delta, L)$--model of a function for optimization tasks, whereby solving a primal problem there is a necessity to recover a solution of a dual problem. The conception of $(\delta, L)$--model is based on the conception of $(\delta, L)$--oracle which was proposed by Devolder–Glineur–Nesterov, herewith the authors proposed approximate a function with an upper bound using a convex quadratic function with some additive noise $\delta$. They managed to get convex quadratic upper bounds with noise even for nonsmooth functions. The conception of $(\delta, L)$--model continues this idea by using instead of a convex quadratic function a more complex convex function in an upper bound. Possibility to recover the solution of a dual problem gives great benefits in different problems, for instance, in some cases, it is faster to find a solution in a primal problem than in a dual problem. Note that primal-dual methods are well studied, but usually each class of optimization problems has its own primal-dual method. Our goal is to develop a method which can find solutions in different classes of optimization problems. This is realized through the use of the conception of $(\delta, L)$--model and adaptive structure of our methods. Thereby, we developed primal-dual adaptive gradient method and fast gradi\-ent method with $(\delta, L)$--model and proved convergence rates of the methods, moreover, for some classes of optimization problems the rates are optimal. The main idea is the following: we find a dual solution to an approximation of a primal problem using the conception of $(\delta, L)$--model. It is much easier to find a solution to an approximated problem, however, we have to do it in each step of our method, thereby the principle of “divide and conquer” is realized.

\includepdf[pages=-]{2020_02_02.pdf}

\section{Статья ``Accelerated and nonaccelerated stocha\-stic gradient descent with model conception''}

\hspace{1cm}

\textbf{Авторы:} Dvinskikh D., Tyurin A., Gasnikov A., Omelchenko S.

\textbf{Аннотация:} In this paper, we describe a new way to get convergence rates for optimal methods in smooth (strongly) convex optimization tasks. Our approach is based on results for tasks where gradients have nonrandom small noises. Unlike previous results, we obtain convergence rates with model conception.

\textbf{Примечание:} статья принята в печать в журнал Mathematical Notes V. 108, № 4. Ниже представлен препринт arXiv: 2001.03443.

\includepdf[pages=-]{2001.03443.pdf}

\section{Статья ``Heuristic adaptive fast gradient method in stochastic opti\-mization tasks''}

\hspace{1cm}

\textbf{Авторы:} Ogaltsov A., Tyurin A.

\textbf{Аннотация:} A heuristic adaptive fast stochastic gradient descent method is proposed. It is shown that this algorithm has a higher convergence rate in practical problems than currently popular optimization methods. Furthermore, a justification of this method is given, and difficulties that prevent obtaining optimal estimates for the proposed algorithm are described.

\textbf{Примечание:} ниже представлена русская версия статьи в Журнал вычислительной математики и математической физики, 2020, том 60, № 7, с. 1143–1150.

\includepdf[pages=-]{VYC1143.pdf}

\section{Статья ``A stable alternative to Sinkhorn’s algo\-rithm for regularized optimal transport''}

\hspace{1cm}

\textbf{Авторы:} Dvurechensky P., Gasnikov A., Omelchenko A., Tyurin A.

\textbf{Аннотация:} In this paper, we are motivated by two important applica\-tions: entropy-regularized optimal transport problem and road or IP traffic demand matrix estimation by entropy model. Both of them include solving a special type of optimization problem with linear equality constraints and objective given as a sum of an entropy regularizer and a linear function. It is known that the state-of-the-art solvers for this problem, which are based on Sinkhorn’s method (also known as RSA or balancing method), can fail to work, when the entropy-regularization parameter is small. We consider the above optimization problem as a particular instance of a general strongly convex optimization problem with linear constraints. We propose a new algorithm to solve this general class of problems. Our approach is based on the transition to the dual problem. First, we introduce a new accelerated gradient method with adaptive choice of gradient’s Lipschitz constant. Then, we apply this method to the dual problem and show, how to reconstruct an approximate solution to the primal problem with provable convergence rate. We prove the rate $O(1/k^2)$, k being the iteration counter, both for the absolute value of the primal objective residual and constraints infeasibility. Our method has similar to Sinkhorn’s method complexity of each iteration, but is faster and more stable numerically, when the regularization parameter is small. We illustrate the advantage of our method by numerical experiments for the two mentioned applications. We show that there exists a threshold, such that, when the regularization parameter is smaller than this threshold, our method outper\-forms the Sinkhorn’s method in terms of computation time.

% \includepdf[pages=-]{CMMP1085.pdf}

\section{Статья ``Dual approaches to the minimization of strongly convex functionals with a simple struc\-ture under affine con\-stra\-ints''}

\hspace{1cm}

\textbf{Авторы:} Anikin A.,  Gasnikov A., Dvurechensky P., Tyurin A., Chernov A.

\textbf{Аннотация:} A strongly convex function of simple structure (for example, separable) is minimized under affine constraints. A dual problem is construct\-ed and solved by applying a fast gradient method. The necessary properties of this method are established relying on which, under rather general conditions, the solution of the primal problem can be recovered with the same accuracy as the dual solution from the sequence generated by this method in the dual space of the problem. Although this approach seems natural, some previously unpublished rather subtle results necessary for its rigorous and complete theoretical substantiation in the required generality are presented.

\includepdf[pages=-]{Anikin et al_2017_Dual Approaches to the Minimization of Strongly Convex Functionals with a Simple Structure under Affine Constraints.pdf}

\section{Статья ``Accelerated gradient sliding for minimizing the sum of functions''}

\hspace{1cm}

\textbf{Авторы:} Dvinskikh D., Omelchenko A., Gasnikov A., Tyurin A.

\textbf{Аннотация:} In this article, we propose a new way to justify the accelerated gradient sliding of G. Lan, which allows one to extend the sliding technique to a combination of an accelerated gradient method with an accelerated variance reduced method. We obtain new optimal estimates for solving the problem of minimizing a sum of smooth strongly convex functions with a smooth regularizer. 

% \includepdf[pages=-]{CMMP1085.pdf}

}  % thesis

\end{document}